\magnification=1200
\centerline{POST-HUMAN MATHEMATICS.}
\bigskip\bigskip
\centerline{by David Ruelle\footnote{$\dagger$}{Math. Dept., Rutgers University, and 
IHES, 91440 Bures sur Yvette, France. email: ruelle@ihes.fr}.}
\bigskip\bigskip
{\sl Abstract:} 
Present day mathematics is a human construct, where computers are used more and more but do not play a creative role.  This situation may change however: computers may become creative, and since they function very differently from the human brain they may produce a very different sort of mathematics.  We discuss what this post-human mathematics may look like, and the philosophical consequences that this may entail.
\vfill\eject
\noindent
{\bf1. Introduction: vitalism.}
\medskip
	The purpose of this note is to look into what may be the mathematics of tomorrow, if and when computer programs have gained a mathematical creativity that they are currently lacking.
\medskip	
	Before that however it will be illuminating to discuss briefly an aspect of the biology and chemistry of yesterday, namely vitalism.  Vitalism is the notion that living matter contains a vital principle which is absent from non-living entities, so that living matter obeys different laws from those that rule non-living matter.  This is an old idea, and it is by no means ridiculous.  This idea has led in chemistry to a distinction between organic and inorganic substances.  The great Swedish chemist J.J. Berzelius (1779-1848) suggested that organic chemical compounds contained a vital force, absent from inorganic compounds.  Indeed, it did not appear possible to synthesize  organic compounds from inorganic ones.  That was until the synthesis of urea by F. W\"ohler in 1828, followed by the synthesis of acetic acid by A. Kolbe in 1845, and the synthesis of more and more organic compounds after that.
\medskip	
	Organic synthesis in the lab has different features from organic synthesis {\it in vivo}, but there appears to be no limit to what can be synthesized.  Even the possibility of creating life in the lab is no longer a taboo subject; many scientists believe that it will be done, although what kind of life and when this will be achieved remain debated questions.
\medskip\noindent
{\bf 2. The uniqueness of human intelligence.}
\medskip
	Although vitalism has largely been abandoned by scientists, there is a belief related to vitalism which remains quite popular, this concerns the human mind and human intelligence.  A lot of people, and this includes very respected scientists, believe that the human mind has unique creative abilities which neither animals nor computers can duplicate.  We shall leave animals out of the present discussion, and concentrate on one aspect of creativity: mathematical creativity.  I see two significant arguments in favor of the uniqueness of human mathematical creativity.  The first argument is that we have an introspective feeling of being able to think (in the present case to think about mathematical problems) which we would deny computers because they are just machines.  The second argument is the lack of serious examples of mathematical creativity by computers.  Let me now discuss and question these two arguments.
\medskip
	{\bf 2.1 Can a computer think?}
\medskip
	When I speak of a computer in what follows, I always mean hardware plus software, i.e., a computer running a suitable program.  The feeling that we can think while computers cannot has been considered by Alan Turing [7].  He discussed the possibility of a test which a human would pass, but which a machine would fail.  The consensus at this time is that no such test can be devised.  A discussion of intelligence based on the introspective feeling that only we can think is thus of limited interest.  We can however say that if present or future computers think, their thinking must be of a very different nature from ours: this question can be objectively analyzed, and has been discussed by Johann von Neumann in a remarkable little book called {\it The computer and the Brain} [3].  Both the computer and the brain are information processing systems, and a detailed comparison made by von Neumann shows that they function very differently.  The brain is slow, prone to error, has limited memory, and is very highly parallel (it is a connected many-channel system).  By comparison, computers are very fast, reliable, have large memory, and are usually not very highly parallel.  Computer intelligence, if it exists, is thus expected to be quite different from human intelligence.
\medskip
	Both computers and humans have several types of memories.  For the task of putting a string of digits in short term memory, a computer is far superior to a human (who can remember only about seven digits).  But our brain also preserves a lifetime of long term memories: an apparently limitless amount.  By comparison, early computers didn't do that well, but things have changed: progress in translation by computers shows that they are progressively becoming able to master the huge corpus of data corresponding to the knowledge of a human natural language.
\medskip
	In everyday life we now frequently have to interact with computer programs, and this influences the opinion we may have on their intelligence.  When I use Google, I have the feeling of a vast intelligence somehow present there.  (This is because Google has rapid and intelligent access to a huge amount of data.)  By contrast, when I recently tried to submit a paper for publication in a couple of scientific journals, and had to fight over the Internet with their editorial programs, I was convinced that I was facing something viciously {\it un}-intelligent.  You know how it is: the thing demands that you  {\tt please enter digital certificate of virginity of grandmother}, or some such nonsense.  Whatever you answer is, the stupid machine replies {\tt certificate invalid please try again}.  Let  us try to overcome the irritation produced by such experiences, and proceed with a serene discussion.
\medskip
	Our conclusion for the moment will be this: we cannot exclude that computers think, but if and when they do it will probably be in a manner very different from that of humans.  The situation is a bit similar to that of organic chemical synthesis: artificial synthesis of organic compounds is not impossible, far from that, but it is generally achieved quite differently from chemical synthesis in living organisms.
\medskip
	{\bf 2.2 Is mathematical creativity by computers possible?}
\medskip
	We come now to the second point of our discussion: the computer's current lack of mathematical creativity.  Computers have actually gained quite a bit of importance in mathematics.  Just think that the great Bernhard Riemann tested some mathematical ideas by long numerical calculations; his modern colleagues often work in a similar way but do their calculations by computer rather than by hand.  Computers are also used in an essential way to provide parts of rigorous proofs: they perform heavy logical or numerical tasks which are beyond human capabilities.  (An example here is the proof of the four color theorem by Kenneth Appel and Wolfgang Haken [1]).  Let me also mention that some definite mathematical creativity has been obtained using Wilf-Zeilberger pairs [4] to produce new identities involving hypergeometric functions.
\medskip
	However, at this time, the closest that computers have come to really doing mathematics is in computer-verified proofs (so-called {\it formal proofs}).  In brief, a human mathematician transforms a human proof of a theorem (like the {\it prime number theorem}) into a sequence of lemmas in a formal language, and the (nontrivial) proof that the lemmas are correct is left to a computer.  For details we refer to a Hales [2] and further papers in a special issue (December 2008) of the AMS Notices (vol. {\bf 55}) on the subject.  Let me just make a few remarks:
\medskip
	(a) nontrivial theorems like the prime number theorem now have a computer-verified proof,
\medskip
	(b) computer-verified formal proofs are a lot more reliable than human proofs, which always contain somewhat imprecise formulations, and sometimes big mistakes (this has become a major concern with very long modern proofs),
\medskip
	(c) part of the computer-verified formal proofs (the proof of lemmas by a computer program) now escapes intuitive grasp by the human mind: these proofs are no longer completely human,
\medskip
	(d) nevertheless, the creative role of computers in computer-verified proofs is minimal, being limited to a combinatorial search for proofs of lemmas along lines programmed by humans.
\medskip\noindent
{\bf 3. What is mathematical creativity?}
\medskip
	This is not the place for a poetic discourse on creativity in general.  Rather, I want to see what creativity implies in the case of mathematics, how it is implemented by humans, and how it might have a non-human implementation.
\medskip
	It is convenient to assume that some basis of mathematics has been accepted: logical rules of deduction and basic axioms.  The axioms may be the Zermelo-Fraenkel-Choice axioms of set theory, or something similar as implemented in a program for computer-verified proofs.  In brief we assume that a common basis for mathematics is accepted by humans and computers.  Doing mathematics is then finding and proving theorems on the basis of the axioms, using accepted rules of logic.
\medskip
	There is a general limitation in doing mathematics that applies both to humans and to computers: a theorem with a short formulation may have an extremely long proof.  This fact, noted by G\"odel, is of logical origin, and related to the incompleteness theorem.
\medskip
	As we have seen, the human brain has limited memory and is prone to error.  A human mathematical text is thus composed of small units (a few lines, where a computer could handle $10^5$ pages).  A great help in obtaining small units is the use of {\it definitions} (for instance the definition of a compact group, or that of complex numbers) which are agreed upon before making a mathematical statement.
\medskip
	Let  me say this again: the human way of doing mathematics is to write a mathematical text, consisting of short pieces which may be definitions or theorems.  Typically there is a main theorem with a long proof, the long proof consisting of definitions and lemmas (the lemmas are little theorems which follow easily from what is already known).
\medskip
	Mathematicians like Hadamard and Poincar\'e have noted that doing mathematics is a combinatorial task: putting pieces together to obtain an interesting theorem.  There are many choices involved in conjecturing an {\it interesting} theorem and putting together the pieces of a proof.  To guide these choices we have a background of results in the published literature, and of theoretical ideas which may be more or less vague or precise.  The published literature increases with time, and the background of theoretical ideas (which define what is an interesting theorem) also changes.  For instance mathematicians are guided by ideas on the natural structures of mathematics; such ideas have been formulated precisely by Bourbaki, or have been later embodied in the theory of categories and functors.  Structural ideas now play a major role in certain areas of mathematics, where some questions will appear natural for structural reasons, and be systematically asked and investigated.  Less precise theoretical ideas consist of analogies, such as the analogy between the theory of C$^*$-algebras and the theory of compact spaces (this analogy comes from the fact that an abelian C$^*$-algebra is precisely the algebra of complex continuous functions on a compact space -- but this does not say exactly how to perform the generalization from the abelian to the non-abelian situation).
\medskip
	To summarize, doing mathematics may be viewed as a succession of guesses and routine verifications.  The guesses are guided by theoretical ideas that evolve with time.  For a more detailed discussion of these questions see my book {\it The Mathematician's Brain} [5].
\medskip
	The above was a description of human mathematics.  In the case of computer-verified mathematics (formal proofs) a part of the guesses is human, but the routine verifications are made by computer, and this involves a nontrivial combinatorial part, i.e., making many low level guesses.  What remains of human creativity are the many higher level guesses, based on theoretical ideas which are not easily formalized to permit their systematic use.
\medskip\noindent
{\bf 4. Limits to mathematical creativity.}
\medskip
	The ability to do mathematics is a recent development in the evolution of the human brain.  Mathematical ability is related (among other things) to the acquisition of language, which is poorly understood.  The ability to speak was clearly favored by evolution, and the same might be said of the ability to count from 1 to 10.  But the ability to do higher mathematics (like studying Galois theory) is another matter, and most people get along successfully without this ability.  One thing that strikes me is the great disparity of performances of the best mathematicians: if one tries to assess quantitatively the contributions of Riemann, G\"odel, or Grothendieck to mathematics, one could say that it is ten to a hundred times greater than that of a ``normal" high-level mathematician.  In other words, the contribution of one of the ``great" mathematicians mentioned above is worth as much as the contribution of ten to a hundred members of the mathematical section of a good academy of sciences (say French, or US).  I think that most mathematical colleagues would agree with this quantitative estimate (although perhaps without enthusiasm).  This is quite different from the situation for 100 m running where the performances of the best racers are quite similar.  To understand the difference, one may appeal to natural selection, which is clearly not the same for running and for doing mathematics, but natural selection arguments are tricky\footnote{*}{I am indebted to Henri Korn for a discussion of this matter.}, and we shall not go further in this direction.
\medskip
	Note, by the way, that a great mathematician is one who does something new, not one who is good at doing again things that have been done before.  We expect therefore that great mathematicians are quite different from each other, so that they can tackle problems in different manners.  This is indeed the case; for example while one slim volume contains the complete works of Riemann, the work of Grothendieck covers many thousands of pages.  Any way one looks at things, the greatest mathematicians are thus very dissimilar people\footnote{*}{The intellectual diversity of people is of course not limited to mathematics.  A glimpse into nonstandard intellectual abilities is provided by some autists, see [7].}: they are not clustered against some natural limit of what can be done by humans in mathematics.  We have seen earlier that there are limits to human creativity imposed by logic and by the structure of the human brain, but now it also appears that individual mathematicians are not close to a universal limit.  The difficulty to put a limit on human mathematical performance suggests that it will also be difficult to put a limit to computer mathematical performance, once computers start to be creative.
\medskip\noindent
{\bf 5. Post-human mathematics.}
\medskip
	I have pointed out that the intellectual ability to do mathematics is a recent development from the point of view of the biological evolution of the human brain.  I find it hard to believe that this recent development has produced something so unique that it cannot be successfully imitated by computers.  I think that the situation of mathematical creativity today  is like that of organic synthesis before W\"ohler: it is one stage in an evolution, and there will be later stages.  The big question is then: what kind of mathematics could be produced by artificial mathematical creativity?  What if we have a computer with access to some sort of mathematical literature, the ability to perform routine proofs, but also the ability to make intelligent guesses because it would have been taught a proper background of theoretical ideas?  What if it developed its own non-human background of theoretical ideas?
\medskip
	Let me interrupt myself here for a brief psychological digression.  I am not eager to see computers replace human mathematicians.  It would, or will, be a sad thing in some respects, but I don't want to shut my eyes to the possibility.  Think of the enormous effects that industrial organic synthesis has had on mankind, some have been nice and some terrible, but there is certainly no way back to the pre-W\"ohler times.  Similarly, mathematics is probably entering soon a completely new era, and we might as well try to guess where this will lead us.
\medskip
	If we assume that a computer has been taught to be mathematically creative, we can then imagine that it could beat human mathematicians at their own game.  This means that it could give proofs of conjectures, or interesting new theorems, which human mathematicians could understand and perhaps admire.  But it is more likely that, once a computer becomes creative, it will do things quite differently from humans.  Here are a couple of possibilities:
\medskip
	(a) The computer could prove an interesting result, but with a proof impenetrable to humans, because it would use long development in some formal language with no reasonably brief translation into familiar human language.  (The Appel-Haken proof of the four color theorem, or the computer verifications using formal proofs, are examples of this).
\medskip
	(b)  The computer could prove an important theorem, but with a statement impenetrable to humans (again because it would have no reasonably brief translation to human mathematical language).  The computer might convince us that this theorem is important, because it implies a number of interesting conjectures that we can understand.  But our brain could not make sense of the theorem itself.
\medskip
	The above possibilities raise big questions on the nature of mathematics.  We can see today's mathematics as a reasonably well structured landscape, with big domains like algebraic geometry or smooth dynamics, and important theorems like the prime number theorem.  This structured mathematical landscape is related to the possibilities of the human brain.  Is there a structure to mathematics which is independent of the human brain?  Could an intelligent computer develop a new mathematical landscape similar to the one we know?  To discuss these questions, remember the logical fact that a theorem with a short formulation may have an extremely long proof.  To develop mathematical knowledge in an economical way, one avoids repeating similar very long proofs.  One tries instead to obtain new results by a relatively short proof using already known results.  The human way of developing mathematics produces thus a network of results related by ``understandable'' proofs (not too long, but also not ``inhumanly'' formal).  The network is constantly being reworked and improved by human mathematicians to ``reveal natural underlying structures''.  This is how the structural landscape that we know for human mathematics has been obtained.  As we have pointed out, there is a logical reason behind this structural landscape, but there are also specifically human reasons (not analyzed in detail but, in brief: the human brain favors short formulations, ``understandable'', and ``interesting'' arguments).  Do we believe that logical factors prevail over human specificities in producing the sort of structural mathematical landscape that we know?
\medskip
	I fear that we must consider another possibility: perhaps computers will develop mathematical abilities so that they can answer efficiently questions that we ask them, but perhaps their efficient way of thinking will have no structural basis recognizable by humans.  If this happens, the superiority of our human intelligence will be strongly challenged: we shall watch an intelligent computer doing mathematics in much the same way as a chimpanzee could watch a human scientist reading a book on Galois theory.\footnote{*}{In this respect, Jean-Pierre Eckmann reminded me of Fred Hoyle's novel {\it The Black Hole} [], which describes human contact with a superior intelligence.}
\medskip
	The gut feeling of many lovers of mathematics will be that they can't believe my chimpanzee story.  Just as Berzelius could not believe that organic compounds would ever be synthesized in the lab.  To go beyond such visceral reactions, I would like to compare mathematics and music.  Mathematicians often love music: there is harmony, beauty, and a sense of infinity in both mathematics and music.  There is also the use of fractions to describe musical intervals (approximately) but this is a somewhat limited relation.  In fact, mathematics and music are conceptually very different things, but it is important for our discussion that they evoke similar esthetic reactions, probably because they involve related activities of the brain.  Here we must remember that there are two sides to mathematics: one is non-human logical necessity, the other is human brain activity.  The non-human logic has little to do with music, but it could be accessible to computers.  As a human brain activity, mathematics is related to other brain activities, and has apparently a privileged relation with music.  The interplay between the human and non-human sides of mathematics is being modified by the irruption of computers into the game.  How this modified interplay will develop in the years to come will be fascinating to observe.
\medskip
	Protagoras has said that ``man is the measure of all things: of things which are, that they are, and of things which are not, that they are not''.  This remains true for us in the sense that everything we know is known to us through our own human brain.  This remains true even though we understand today that the planet Earth on which we live is an infinitesimal speck of dust in the physical universe.  This will remain true tomorrow even if we find that human mathematical discoveries are dwarfed and humbled by computer mathematics.
\medskip
	The present text is a slightly reworked version of a presentation made at the 20-th anniversary conference of the ESI in Vienna, 29 April 2013
\medskip\noindent
{\bf References.}
\medskip
	[1] K. Appel and W. Haken.  {\it Every Planar Map is Four-Colorable.}  AMS, Providence, 1989.
\medskip
	[2] T.C. Hales.  ``Formal proof.''  AMS Notices {\bf 55},1370-1380(2008).
\medskip
	[] F. Hoyle.  {\it The Black Cloud.}  William Heinemann, London, 1957.
\medskip
	[3] J. von Neumann.  {\it The Computer and the Brain.}  Yale U. P., New Haven, 1958.
\medskip
	[4] M. Petkovsek, H. Wilf and D. Zeilberger.  {\it A+B.}  A K Peters, 1996.
\medskip
	[5] D. Ruelle.  {\it The Mathematician's Brain.}  Princeton U. P., Princeton, 2007.
\medskip
	[6] D. Tammet.  {\it Born on a Blue Day.}  Hodder and Stoughton, London, 2006.
\medskip
	[7] A. Turing.  ``Computing machinery and intelligence.''  Mind {\bf 59},433-460(1950).

\end